\documentclass[11pt]{amsart}
\usepackage{amsmath, amsthm, amsfonts, ifpdf}
\usepackage[dvipsnames,usenames]{color}
\usepackage[demo]{graphicx}
\usepackage{multicol}
\usepackage{tikz}
\usepackage{pgfmath}
\usetikzlibrary{arrows,shapes,trees,backgrounds}
\usetikzlibrary{intersections}

\theoremstyle{plain}
\newtheorem{theorem}{Theorem}[section]
\newtheorem*{theorem*}{Theorem}

\newtheorem{pro}[theorem]{Proposition}
\newtheorem{Def}[theorem]{Definition}
\newtheorem{lem}[theorem]{Lemma}

\theoremstyle{definition}
\newtheorem*{Def*}{Definition}

\newtheorem{Rem}[theorem]{Remark}

\numberwithin{equation}{section}

\newcommand{\bpo}{\begin{pro}}
\newcommand{\epo}{\end{pro}}
\newcommand{\be}{\begin{equation}}
\newcommand{\ene}{\end{equation}}
\newcommand{\br}{\begin{Rem}}
\newcommand{\er}{\end{Rem}}
\newcommand{\bl}{\begin{lem}}
\newcommand{\el}{\end{lem}}
\newcommand{\bd}{\begin{Def}}
\newcommand{\ed}{\end{Def}}
\newcommand{\ben}{\begin{enumerate}}
\newcommand{\een}{\end{enumerate}}
\newcommand{\bp}{\begin{proof}}
\newcommand{\ep}{\end{proof}}
\newcommand{\beq}{\begin{equation*}}
\newcommand{\eeq}{\end{equation*}}
\newcommand{\bear}{\begin{eqnarray*}}
\newcommand{\eear}{\end{eqnarray*}}
\newcommand{\bt}{\begin{theorem}}
\newcommand{\et}{\end{theorem}}
\newcommand{\bst}{\begin{split}}
\newcommand{\est}{\end{split}}

\newcommand{\bal}{\begin{aligned}}
\newcommand{\eal}{\end{aligned}}

\renewcommand{\P}{\partial}
\newcommand{\F}[2]{\frac{#1}{#2}}
\newcommand{\la}{\langle}
\newcommand{\ra}{\rangle}

\newcommand{\R}{\mathbb{R}}

\newcommand{\bnb}{\bar{\nabla}}
\newcommand{\nb}{\nabla}

\newcommand{\ta}{\theta}

\renewcommand{\H}{\mathbb{H}}
\newcommand{\Ta}{\Theta}

\newcommand{\PLH}{{\mkern-1mu\times\mkern-1mu}}

\newcommand{\NR}{N^2\PLH\R}
\newcommand{\TNR}{\widetilde{N^2\PLH\R}}
\newcommand{\lh}{\tilde{h}}
\newcommand{\tM}{\tilde{M}}
\begin{document}
\title[Complete graphs in $\NR$]{ The boundary behavior of domains with complete translating, minimal and CMC graphs in $N^2\PLH \R$}
\author{Hengyu Zhou}
\address{Department of Mathematics, Sun Yat-sen University, No. 135, Xingang Xi Road, Guangzhou, 510275, People's Repulic of China}
\email{zhouhy28@mail.sysu.edu.cn}
\date{\today}
\begin{abstract} 
In this note we discuss graphs over a domain $\Omega\subset N^2$ in the product manifold $N^2\PLH\R$. Here $N^2$ is a complete Riemannian surface and $\Omega$ has piece-wise smooth boundary. Let $\gamma \subset \P\Omega$ be a smooth connected arc and $\Sigma$ be a complete graph in $N^2\PLH\R$ over $\Omega$.  We show that if $\Sigma$ is a minimal or translating graph, then $\gamma$ is a geodesic  in $N^2$. Moreover if $\Sigma$ is a CMC graph, then $\gamma$ has constant principle curvature in $N^2$. This explains the infinity value boundary condition upon domains having Jenkins-Serrin theorems on minimal and CMC graphs in $N^2\PLH\R$. 
\end{abstract}
\maketitle
\section{Introduction}
In this paper we are interested in the asymptotic behavior of translating, minimal and constant mean curvature (CMC) graphs over a domain in a Riemannian surface. The purpose is to
establish the connection between the completeness of those graphs over a domain and the property of its boundary. We are motivated by recent progresses on complete translating graphs in $\R^3$ and the Jenkins-Serrin theory on minimal graphs and CMC graphs.\\
\indent Before giving more details let us introduce the concept of translating graphs. We apply the following notation throughout this paper: $N^2$ is a complete Riemannian surface with a metric $\sigma$, $N^2\PLH\R$ is the product manifold $\{(x,r): x\in N^2, r\in \R\}$ equipped with the metric $\sigma+dr^2$ and $\Omega$ is a domain in $N^2$ with piecewise smooth boundary.\\
\indent A \textit{translating graph} in $N^2\PLH\R$ if it is the graph of $u(x)$ where $u(x):\Omega\rightarrow \R$ is the solution of a mean curvature type equation given as follows: 
\be \label{def:tsgraph}
div(\F{Du}{\sqrt{1+|Du|^2}})=\F{1}{\sqrt{1+|Du|^2}}
\ene 
where $Du$ is the gradient of $u$ and div is the divergence of $N^2$. Translating surfaces characterize the type II finite singularity of mean curvature flow in Euclidean space (see \cite{AngV97}, \cite{Ang95} and \cite{HS08}). Some geometric  properties were investigated in \cite{AW94, WXJ11, GJJ10, CSS07, Sj16} etc.\\
\indent  Recently Shahriyari \cite{Sha15} showed that if $\Sigma$ is a complete translating graph over a smooth domain $\Omega\subset \R^2$ in $\R^3$, then $\P\Omega$ has to be a geodesic. This raises a question as follows.\\
\indent \emph{ Is there a connection between the completeness of graphs with certain properties over a domain and the boundary behavior of this domain? }\\
\indent An answer for this question in the case of CMC graphs was already founded by Spruck (section 8 in \cite{spr72}). One of his results says that suppose the function $u(x)$ in a domain $\Omega$ in $\R^n$ goes to $+\infty$ uniformly as $x$ approaches to an connected open domain $\Gamma\subset\P\Omega$ and the graph of $u(x)$ is a (complete obviously) CMC graph in $\R^{n+1}$, then $\Gamma$ has constant mean curvature in $\R^n$. \\
\indent Our question is also related to the Jenkins-Serrin theory on minimal graphs and CMC graphs in product manifolds (see Jekins-Serrin\cite{JS68}). For an excellent summary of this topic we refer to Eichmair-Metzger \cite{EM16}. Its basic setting is given as follows. Let $\Omega$ be a domain in $N^2$ with its boundary $\P\Omega$ which is composed with $\P_{+}\Omega$, $\P_{-}\Omega$ and $\P_{0}
\Omega$. The Jenkins-Serrin theory seek to a smooth function $u(x)$ on $\Omega$ such that $\Sigma$, the graph of $u(x)$,  is minimal or of CMC in $N^2\PLH\R$ and $u(x)$ approaches to $+\infty(-\infty)$ when $x$ is close to $\P_{+}\Omega (\P_{-}\Omega)$ and approaches to continuous data when $x$ is close to $\P_{0}\Omega$. Generally this theory also requires that $\P_{+}\Omega$ and  $\P_{-}\Omega$ are minimal or have constant principle curvature in $N^2$ respectively (see \cite{spr72},\cite{PL09} and \cite{EM16}). One interesting application of Jenkins-Serrin theorems is the construction of a harmonic diffeomorphism from the complex plane $\mathbb{C}^2$ to the hyperbolic plane $\H^2$ by Collin-Rosenberg \cite{CR10}. \\
\indent Now our main result will answer the question mentioned above. It also explains the conditions on $\P_{+}\Omega$ and $\P_{-}\Omega$ in the Jekins-Serrin theory. We say that a graph $\Sigma$ in $N^2\PLH \R$ over $\Omega$ is complete approaching to a connected arc $\gamma\subset \P\Omega$ if $\Sigma$ can not be extended along $\gamma$ as a complete graph over a neighborhood of  $\gamma$. The main result of this paper is stated as follows. 
\bt\label{thm:MT1}(\emph{Theorem \ref{thm:geo}})
Let $N^2$ be a complete Riemannian surface and $\Omega\subset N^2$ is a domain with piecewise smooth boundary. Let $\gamma \subset \P\Omega$ denote a smooth connected arc and let $\Sigma$ be the graph of a smooth function $u(x)$ on $\Omega$ in the product manifold $N^2\PLH\R$.  \\
\indent  Suppose $\Sigma$ is complete approaching to $\gamma$. Then we have 
\begin{enumerate} 
	\item if $\Sigma$ is a translating or minimal graph, then $\gamma$ is a geodesic arc;
	\item if $\Sigma$ is a CMC graph, then $\gamma$ has constant principle curvature. 
\end{enumerate} Moreover only one of the following holds: (1) $u(x)\rightarrow +\infty$ as $x\rightarrow x_0$ for all $x_0\in \gamma$; (2) $u(x)\rightarrow -\infty$ as $x\rightarrow x_0$ for all $x_0\in \gamma$. 
\et
\br 
The reason that we only work in a surface $N^2$ is that we need curvature estimates of stable type surfaces in three manifolds (see Section 4). We are working on a project that deals with the higher dimension version of Theorem \ref{thm:MT1}.  \er  
The essential part in the proof of our main result is when $\Sigma$ is a translating graph. The other two cases can be achieved with  minor modification (see Section 5). Its basic idea is inspired from Shariyari \cite{Sha15}. \\
\indent When $\Sigma$ is a translating graph, we show that $\Sigma$ is stable and minimal with respect to a weighted product metric (see Theorem \ref{thm:mta}). The curvature estimate of stable minimal surfaces in three dimensional manifolds ( Schoen \cite{Soe83} and Minicozzi-Colding \cite{CM02}) gives a family of simply connected disks on $\Sigma$ with fixed diameter $\delta$ (see Lemma \ref{lm:ti}) centered at points $(x_n, u(x_n))$ where $x_n$ goes to a point in $\gamma$. These disks has a vertical limit $F$ according to Theorem \ref{thm:est2}. Moreover $F$ is minimal since the completeness of $\Sigma$ guarantees that the angle function on $F$ has to vanish according to Theorem \ref{thm:convergence} (see Lemma \ref{lm:minimal}). This implies that $\gamma$ is a geodesic. Notice that when $\Sigma$ is a CMC graph, the curvature estimate we need is from Zhang \cite{Zhang05} (see Theorem \ref{thm:cmc:ce}).\\
\indent Our paper is organized as follows. In Section 2 we show that a translating graph is minimal and stable with respect to a weighted metric in $N^2\PLH\R$. In Section 3 we compute the sectional curvature of this weighted metric. In Section 4  all curvature estimates of stable minimal surfaces and CMC surface that we need are collected. In Section 5 we prove Theorem \ref{thm:MT1}. In appendix A we construct translating graphs in $N^2\PLH\R$ where $N^2$ has certain warped product structure. A particular example is that $N^2$ is the two dimensional hyperbolic space $\H^2$. \\

\section{Stability}
Let $\Sigma$ be a translating graph of $u(x)$ in $\NR$ where $u(x)$  satisfies
\eqref{def:tsgraph} on a domain $\Omega$. We follow the notation in \cite{HZ16}. The upward normal vector $\vec{v}$ is $\Ta(\P_r-Du)$ where $Du$ is the gradient of $u(x)$ with respect to $N^2$. Suppose $\{\P_1,\P_2\}$ is a local frame on $N^2$. We denote $\P_i+u_i\P_r$ by $X_i$ for $i=1,2$. Then $\{X_1, X_2\}$ is a local frame on $\Sigma$. The second fundamental form of $\Sigma$ is 
$$
h_{ij}=-\la \bnb_{X_i}X_j,\vec{v}\ra
$$
With these notation one sees that 
$
H=-div(\F{Du}{\sqrt{1+|Du|^2}})
$. Therefore an equivalent form of \eqref{def:tsgraph} is
\be\label{eq:angle}
H=-\Ta=-\la\vec{v},\P_r\ra
\ene
where $\Ta$ is referred as the angle function of $\Sigma$. The corresponding Codazzi equation and Gauss equation take the following form:
\begin{gather}
R_{ijkl}=\bar{R}_{ijkl}+(h_{ik}h_{jl}-h_{il}h_{jk})\\
h_{ij,k}=h_{ik,j}+\bar{R}_{\vec{v}ijk}
\end{gather} where $R$ and $\bar{R}$ are the Riemann curvature tensor of $\Sigma$ and $N^2\PLH\R$ respectively. \\
\indent 
Let $\TNR$ denote the manifold $\{(x,r): x\in N^2, r\in \R\}$ equipped with a weighted product metric $e^{r}(\sigma+dr^2)$.\\
\indent First we show that
\bt\label{thm:mta} Let $u(x)$ be a solution in \eqref{def:tsgraph}. Its graph $\Sigma=(x,u(x))$ is a stable minimal surface in $\TNR$.
\et
\bp
Before the proof, let us check the area functional of $\TNR$ given by
$$
F(\Sigma)=\int_{\Sigma}e^{r}d\mu
$$
where $d\mu$ is the volume of $\Sigma$ in product manifold $\NR$. Let $\Sigma_s$ be a family of surfaces satisfying
\be
\F{\P\Sigma_s}{\P s}|_{t=0}=\phi\vec{v}\quad \text{with}\quad    \Sigma_0=\Sigma
\ene
where $\phi(x)$ is a smooth function on $\Sigma$ with compact support. We view $\Sigma_s$ as a curvature flow of $\Sigma$ in $\NR$. From the classical computation in curvature flows (see Huisken-Polden \cite{HP99}), we have
\be\label{eq:pd}
\begin{split}
	\F{\P\vec{v}}{\P s}|_{s=0}&=-\nb\phi\\
	\F{\P H}{\P s}|_{s=0}&=-\Delta\phi-(|A|^2+\bar{R}ic(\vec{v},\vec{v}))\phi
\end{split}
\ene
where  $\nb$, $\Delta$ are the covariant derivative and Laplacian of $\Sigma$ in $\NR$, and $\bar{R}ic$ is the Ricci curvature tensor of $\NR$. According to \eqref{eq:pd} and \eqref{eq:angle}, a direct computation shows that
\be
\begin{split}
	\F{\P F(\Sigma_s)}{\P s}|_{s=0} &=\int_{\Sigma}\phi(H+\la \vec{v},\P_r\ra)e^{r} d\mu=0\\
	\F{\P^2 F(\Sigma_s)}{\P^2 s}|_{s=0}&=-\int_{\Sigma}\phi(\Delta\phi+(|A|^2+\bar{R}ic(\vec{v},\vec{v}))\phi+\la \nb \phi,\P_r\ra)e^{r}d\mu
\end{split}
\ene
For convenience of computation, we define an elliptic operator $L$ as follows:
\be\label{def:L}
L\phi=\Delta\phi+(|A|^2+\bar{R}ic(\vec{v},\vec{v}))\phi+\la \nb \phi,\P_r\ra
\ene
With this notation it is sufficient to check that whether 
\be\label{eq:second:variation}
\F{\P^2 F(\Sigma_s)}{\P^2 s}|_{s=0}=-\int_{\Sigma}\phi L\phi e^{r}d\mu
\ene
is negative. Since $\Sigma$ is a graph, its angle function$\Ta=\la \vec{v},\P_r\ra >0$. Thus we can write $\phi=\eta\Ta$ where $\eta$ is another function over $\Sigma$ with compact support. Therefore, we obtain that
\be \label{eq:rst}
\phi L\phi=\eta\Ta(\eta L\Ta+\Ta\Delta \eta+2\la \nb \eta, \nb\Ta\ra +\Ta\la\nb\eta,\P_r\ra)
\ene
The reason we adapt this form is based on a general
formula of  $\Delta\Ta$ as follows. 
\bl \label{lm:graphic}On any $C^2$ surface $S$ in $\NR$, it holds that
\be \label{eq:reta}
\Delta \Ta+(|A|^2+\bar{R}ic(\vec{v},\vec{v}))\Ta-\la \nb H, \P_r\ra =0
\ene
where $A$ is the second fundamental form of $S$.
\el
\bp Fix a point $p\in S$. Choose an orthonormal frame $\{e_1, e_2\}$ on $S$ such that $\nb_{e_i}e_j(p)=0$ and $\la e_i, e_j\ra =\delta_{ij}$.\\
\indent Then $\bnb_{e_i}e_j(p)=-h_{ij}\vec{v}$ where $\bnb$ denotes the covariant derivative of $\NR$ and $\vec{v}$ is the normal vector of $S$. Since $\NR$ is a product manifold, it is well-known that $\bnb_X\P_r=0$ for any smooth  vector field $X$. We compute $\Delta \Ta$ as follows.
\begin{align}
\Delta \Ta(p)&=\nb_{e_i}\nb_{e_i}\la \P_r,\vec{v}\ra-\nb_{\nb_{e_i}e_i}\Ta(p)\notag\\
&=e_i\la \P_r, h_{ik}e_k\ra (p)\notag\\
&=h_{ik,i}\la\P_r,e_k\ra-|A|^2\Ta\label{eq:basic}
\end{align}
Recall that the Codazzi equation (Chapter 6 in \cite{doC92}) says that
\be
h_{ik,i}=h_{ii,k}+\bar{R}(\vec{v},e_i,e_k,e_i)
\ene
where $\bar{R}$ denotes the Riemann curvature tensor of $\NR$. Thus
\be
h_{ik,i}\la\P_r,e_k\ra=\la \nb H,\P_r\ra+\bar{R}ic(\vec{v},\la \P_r, e_k\ra e_k)
\ene
We observe that
\be
\la \P_r, e_k\ra e_k=\P_r-\Ta \vec{v}
\ene
and $\bar{R}ic(\vec{v},\P_r)=0$ because $\bnb_{X}\P_r=0$ for any vector $X$. This implies that
$$
h_{ik,i}\la\P_r,e_k\ra=\la \nb H,\P_r\ra-\bar{R}ic(\vec{v},\vec{v})\Ta
$$ Combining this with \eqref{eq:basic}, we achieve the lemma.
\ep
Now we go back to the proof of Theorem \ref{thm:mta}. By assumption $\Sigma$ is a translating graph in $\NR$, then $H=-\Ta$. Hence \eqref{eq:reta} is written as
\be
L\Ta=0
\ene
and therefore \eqref{eq:rst} becomes that
\be
\phi L\phi =\eta\Ta(\Ta\Delta \eta+2\la \nb \eta, \nb\Ta\ra -\Ta\la\nb\eta,\P_r\ra)
\ene
On the other hand, the divergence of $\eta \Ta^2\nb \eta e^{r}$ is computed as follows.
\be
\begin{split}
	div(\eta\Ta^2\nb \eta e^{r})&=\eta\Ta e^{r}(\Ta\Delta \eta+2\la \nb \eta, \nb\Ta\ra +\Ta\la\nb\eta,\P_r\ra)+\Ta^2|\nb \eta|^2 e^{r}\\
	&=\phi L\phi e^{r} +\Ta^2|\nb \eta|^2 e^{r}
\end{split}
\ene
Combining this expression with \eqref{eq:second:variation} and applying 
the divergence theorem we obtain that 
$$
\F{\P^2 F(\Sigma_s)}{\P^2 s}|_{s=0}=\int_{\Sigma} \Ta^2|\nb \eta|^2 e^{r}d\mu\geq 0
$$
Then we conclude that  $\Sigma$ is stable and minimal in $\TNR$.
\ep
Lemma \ref{lm:graphic} leads to a rigidity result of limit surfaces for the $C^2$ convergence of minimal graphs, translating graphs and CMC graphs in $N^2\PLH\R$. It is an important ingredient in the proof of Lemma \ref{lm:minimal} (see Section 5). A similar result appeared in Lemma 2.3 of Eichmair \cite{Ecm10} in the setting of marginally outer trapped
surfaces. 
\bt \label{thm:convergence} Let $\{\Sigma_n\}_{n=1}^\infty$ be a sequence of smooth connected graphs in $N^2\PLH\R$ with diameter $\delta$ converging uniformly to a connected surface $\Sigma$ in the $C^2$ sense. If all $\Sigma_n$ are minimal in the interior of $\Sigma$ the angle function $\Ta$ satisfies that $\Ta>0$ or $\Ta\equiv 0$. The conclusion is also in the case of minimal or CMC graphs. 
\et 
\bp  Without loss of generality we assume $\Ta>0$ in the interior of all $\Sigma_n$. \\
\indent First we assume $\Sigma_n$ are minimal or CMC. Then $\nb H\equiv 0$. By Lemma \ref{lm:graphic}, we have 
\be
\Delta \Ta+(|A|^2+\bar{R}ic(\vec{v},\vec{v}))\Ta=0
\ene 
on all $\Sigma_n$.
Notice that $\bar{R}ic(\vec{v},\vec{v})\geq -\beta$ here $\beta$ is a positive constant only depending on $N^2$. Then we have $
\Delta \Ta\leq \beta \Ta
$
on all $\Sigma_n$. Since $\Sigma$ are the $C^2$ uniform limit of $\Sigma_n$ as $n\rightarrow \infty$, then it has the property that $\Ta\geq 0$ and satisfies that  $\Delta \Ta\leq \beta \Ta$. By the strong maximum principle of elliptic equations, $\Ta\equiv 0$ or $\Ta>0$ on $\Sigma$.  \\
\indent Assume $\Sigma_n$ are translating graphs. Then $H\equiv -\Ta$ by \eqref{eq:angle}. It also holds that 
$\Delta \Ta\leq \beta \Ta+\la\nb \Ta, \P_r\ra$ on $\Sigma_n$ and $\Sigma$. Based on the strong maximum principle and $\Ta\geq 0$ on $\Sigma$, we obtain the conclusion with a similar derivation as above.
\ep
\section{Sectional Curvature of $\TNR$}
In this section we compute the sectional curvature of $\TNR$ (see Lemma \ref{lm:section}). It should be a classical fact in Riemannian geometry. However we did not find appropriate literatures. For convenience of readers we inlcude its proof here.\\
\indent We first work with the general setting. Suppose $M$ is a Riemannian manifold with a smooth metric $g$. In a local coordinate $\{x_1,\cdots,x_n\}$ , its metric is expressed as
$$
g=g_{ij}dx^idx^j
$$
The Christoffel symbols with this local coordinates are defined by $\nb_{\P_i}\P_j=\Gamma_{ij}^k\P_k$ and computed by the following expressions:
\be
\Gamma_{ij}^k=\F{1}{2}g^{kl}\{\P_ig_{lj}+\P_jg_{li}-\P_l g_{ij}\}
\ene
and the Riemannian curvature tensor is computed by 
\be \label{eq:R}
R(\P_i,\P_j,\P_k,\P_l)=(\P_j\Gamma_{ik}^r-\P_j\Gamma_{ik}^r +\Gamma_{ik}^m\Gamma_{mj}^r-\Gamma_{jk}^m\Gamma_{mi}^r)g_{rl}
\ene
Let $\tM$ denote the same smooth manifold $M$ equipped with a weighted metric $e^{2f}g$ where $f$ is a smooth function on $M$. According the definition above, the corresponding Christoffel symbols $\tilde{\Gamma}_{ij}^k$ satisfy that
\be \label{eq:st}
\tilde{\Gamma}_{ij}^k=\Gamma_{ij}^k+(\delta_{kj}\P_i f+\delta_{ki}\P_j f-g^{kl}\P_lf g_{ij})
\ene
where $\nb$ is the covariant derivative of $M$. \\
\indent We consider the second fundamental form of hypersurfaces in Riemanian manifolds and their conformal deformations. Notice that all computation belows are valid for any dimension. We have the following result.
\bl If $\Sigma$ is a hypersurface in $M$, then it is a hypersurface in $\tM$ and vice versa. Let $(h_{ij})$ and $(\lh_{ij})$ be the second fundamental form of a hypersurface $\Sigma$ in $M$ and $\tM$ respectively. Then they satisfy that
\be\label{eq:second}
\lh_{ij}=e^f(h_{ij}+df(\vec{v})g_{ij})\quad  \lh_{i}^j=e^{-f}(h_i^j+df(\vec{v})\delta_{ij                         })
\ene
where $\vec{v}$ is the normal vector of $\Sigma$ in $M$.
\el
\bp We add $\sim$ for all geometric quantities related on $\tM$. Fix a point $p$ on $\Sigma$. Since $\Sigma$ is a manifold, we can choose a local chart $\{x_1,\cdots, x_n, x_{n+1}\}$ near $p$ in $M$ such that $\{\P_1,\cdots,\P_n\}$ is a local frame on $\Sigma$ and $\vec{v}(p)=\P_{n+1}$. Notice that $e^{-f}\vec{v}$ is the normal vector of $\Sigma$ in $\tM$. According to the definition and applying \eqref{eq:st} one observes that
\begin{align}
\lh_{ij}(p)&=-\tilde{g}(\bnb_{\P_i}\P_{j},e^{-f}\vec{v})=-a_l e^{-f}\tilde{\Gamma}_{ij}^{n+1}\\
&=-a_l e^{f}\Gamma_{ij}^{n+1}+a_l e^{f}
\P_l f g_{ij}\\
&=e^f (h_{ij}+df(\vec{v}) g_{ij})
\end{align}
because $\tilde{g}_{ij}=e^{2f}g_{ij}$.
Applying $\tilde{g}^{ij}=e^{-2f}g^{ij}$, we obtain \eqref{eq:second} from $\lh_i^i=\tilde{g}^{ik}\lh_{ki}$.
\ep
Now we obtain the relationship between the sectional curvature of $\NR$ and this of $\TNR$ by applying the Codazzi equation. 
\bl\label{lm:section} Let $\{\P_1,\P_2,\P_3=\P_r\}$ be a local orthonormal frame on $\NR$. Let $K_{ij}$ and $\tilde{K}_{ij}$ denote the sectional curvature of $\NR$ and $\TNR$ respectively. Then it holds that
\be
\tilde{K}_{ij}(x,r)=e^{-r}(K_{ij}(x)-\F{1}{4})\quad \tilde{K}_{i3}=0
\ene
for $i,j\in\{1,2\}$.
\el
\bp Let $f(r)=\F{r}{2}$. Fix any $r$. Notice that the slice $N^2\PLH \{r\}$ is totally geodesic in $\NR$. Assume $i,j\in\{1,2\}$.  According to \eqref{eq:second}, its second fundamental form in $\TNR$ is
\be
\lh_{ij}=-\F{1}{2}e^{\F{r}{2}}\sigma_{ij}
\ene
It is easy to see that the Riemannian curvature tensor of $N^2\PLH \{r\}$ with respect to the induced metric is $e^{r}R$ where $R$ is the Riemannian curvature tensor of $N^2$. By the Codazzi equation,
$$
\bar{R}_{ijij}(x,r)=e^{r}R_{ijij}(x)-\F{1}{4}e^{r}(\sigma_{ii}\sigma_{jj}-\sigma_{ij})(x)
$$
From $\tilde{g}_{ij}=e^{r}\sigma_{ij}=e^{r}g_{ij}$, a direct computation yields the expression of $\tilde{K}_{ij}$. With a  straightforward computation, we have
$$
\tilde{\Gamma}_{k3}^l=-\F{1}{2}\delta_{kl}
$$
where $\{k,l\}\in \{1,2\}$. According to \eqref{eq:R} we have $\tilde{K}_{i3}=\bar{R}_{i3i3}=0$.
\ep
\section{Curvature Estimates}
We recall some curvature estimates from Schoen \cite{Soe83}, Minicozzi-Colding \cite{CM02} on the stable minimal surfaces $\Sigma$ and Zhang \cite{Zhang05} on stable CMC surfaces immersed in three Riemannian manifolds $M^3$ with sectional curvature $K_M$. For a fixed point $x\in M$, $B_1(x)$ denotes the extrinsic ball in $M$ centered at $x$ with radius $r$. Similarly, For a fixed point $x\in \Sigma$, $B_r^\Sigma(x)$ denotes the intrinsic ball on $\Sigma$ centered at $x$ with radius $r$.
\bt\label{thm:est}(\emph{\cite{Soe83} and \cite{CM02}}) Suppose $\Sigma\subset M^3$ is a stable minimal surface with trivial normal bundle and $B_{r_0}^\Sigma(p)\subset \Sigma\backslash \P\Sigma$ where $|K_M|\leq k^2$ where $r_0<\min\{\F{\pi}{k},k\}$. Then for some positive constant $C=C(k)$ and all $0<r <r_0$,
\be
\sup_{B^\Sigma_{r_0-r}(p)}|A|^2\leq Cr^{-2}
\ene
\et
Now we derive the curvature estimate for translating graphs. The idea follows from  Shariyari \cite{Sha15}. Let $i_0$ denote the injective radius in $N^2$. Without loss of generality, we assume $i_0\leq 1$.
\bt\label{thm:est2} Let $U$ be an open domain in $N^2$ with sectional curvature satisfying $|K_{N^2}(x)|+\F{1}{4}\leq k^2$ for all $x\in U$. Let $\Sigma=(x,u(x))$ be a translating graph in $\NR$ where $x\in U$. If $B^\Sigma_{r_0e^{-1}}(p)\in\Sigma\cap B_1(p) \backslash \P(\Sigma\cap B_1(p))$ and $r_0\leq \min\{\F{\pi}{k\sqrt{e}},k\sqrt{e}, i_0,1\}$ then for some positive constant $C=C(k)$ and all $0<r\leq r_0e^{-1}$,
\be\label{eq:third}
\sup_{B^\Sigma_{r_0e^{-1}-r}(p)}|A|^2\leq Cr^{-2}
\ene
\et
\bp  Fix a point $p=(x_0,y_0)\in \NR$ where $x_0\in U$. Let $B_{r}(p)$ be the ball in $\NR$ containing all points of which the distance to $p$ is $r_0$. Then for any point $(x,y)\in B_{r_0}(p)$ we have
\be\label{eq:relation}
|y-y_0|\leq r_0\leq 1
\ene
\indent  Let $\tilde{B}_{r_0}(p)$ denote the ball $B_{r_0}(p)$ equipped with the conformal metric $e^{r-y_0}(\sigma+dr^2)$. We claim that the sectional curvatures $\tilde{K}$ of all points in $\tilde{B}_{r_0}(p)$ satisfies that
\be\label{eq:sectional}
|\tilde{K}|\leq ek^2
\ene
. By Lemma \ref{lm:section}, the sectional curvature $\tilde{K}_{ij}$ of the ball $B_{r_0}(p)$ equipped with the conformal metric $e^{r}(\sigma+dr^2)$ is $ e^{-r}(K_{ij}-\F{1}{4})$ for $i=1,2$ and $\tilde{K}_{i3}=0$. Multiplying the constant factor $e^{-y_0}$, the sectional curvature of $\tilde{B}_{r_0}(p)$ shall satisfy that
\be\label{eq:sec:est}
\tilde{K}_{ij}(x,r)=e^{y_0-r}(K_{ij}(x)-\F{1}{4})\quad \tilde{K}_{i3}=0;
\ene
for $i,j=1,2$. Combining \eqref{eq:sec:est} with \eqref{eq:relation} and $|K_N(x)|+\F{1}{4}\leq k^2$ yields that  \eqref{eq:sectional}. We obtain the claim.\\
\indent By Theorem \ref{thm:mta}, $\Sigma$ is a stable minimal graph with respect to the metric $e^{r}(\sigma+dr^2)$. Since $y_0$ is a constant,  $\Sigma$ is still a stable minimal graph in $\tilde{B}_{r_0}(p)$ with respect to the metric $e^{r-y_0}(\sigma+dr^2)$. Here multiplying a positive constant $e^{-y_0}$ on the metric does not change the minimal and stable properties of hypersurfaces. Now applying Theorem \ref{thm:est}, the second fundamental form $\tilde{A}$ of $\Sigma$ in $\tilde{B}_{r_0}(p)$ satisfies
\be
\sup_{\tilde{B}_{r_0-\sigma'}(p)}|\tilde{A}|^2\leq C(k)\sigma^{-2}
\ene
where $r_0\leq \min\{\F{\pi}{k\sqrt{e}},k\sqrt{e}\}$, $\sigma'<r_0$ and $C(k)$ is a constant depending on $k$. According to Lemma \ref{eq:second}, we have
\be
|\tilde{A}|^2=e^{-(r-y_0)}|A|^2\geq e^{-1}|A|^2
\ene
where $A$ is the second fundamental form of $\Sigma$ with respect to the metric $\sigma+dr^2$. Similarly, the ball $\tilde{B}_{(r_0-\sigma')}(p)$ with respect to the metric $e^{r-y_0}(\sigma+dr^2)$ contains the ball $B_{(r_0 e^{-1}-\sigma'e^{-1})}(p)$ because of \eqref{eq:relation}.  This yields that
\be
\sup_{B_{(r_0e^{-1}-\sigma'e^{-1})}(p)}|A|^2\leq C(k)e^3(e\sigma')^{-2}
\ene
Let $r$ be $\sigma'e^{-1}$, we obtain \eqref{eq:third}. 
\ep
Following from Zhang \cite{Zhang05}, a CMC surface $\Sigma$ is \emph{stable} if for any $f\in C^\infty_0(\Sigma)$ it holds that $-Lf\geq 0$ where $L$ is given by \eqref{def:L}. Thus a CMC graph in $N^2\PLH \R$ is stable by combining \eqref{def:L} with \eqref{eq:reta}. The curvature estimate of stable CMC surfaces is given as follows. 
\bt\label{thm:cmc:ce}(\emph{Theorem 1.1 in \cite{Zhang05}}) Let $\Sigma$ be an immersed stable CMC $H_0$-surface with trivial normal bundle in a complete three dimensional manifold $M$ where its sectional curvature satisfies $|K_M|\leq k^2$. There exists a positive constant $r_0=r_0(H_0, k, M)$ such that for all $\sigma \in (0, r_0)$ and any $x\in \Sigma$ with geodesic ball $B_{r_0}(x)\bigcap \P\Sigma=\emptyset$ we have for 
$$
\sup_{B^\Sigma_{r_0-\sigma}(x)}|A|^2 \leq C\sigma^{-2}
$$
where $C$ is a constant only depending on $H_0,k$ and $M$.     

\et 
\br In this paper we do not need the precise expressions of those constants.
\er
\section{Proof of the main theorem}
In this section we show the main theorem. 
\bt\label{thm:geo}(\emph{Theorem \ref{thm:MT1}})  Let $N^2$ be a Riemannian surface and $\Omega\subset N^2$ be a domain with piecewise smooth boundaries. Let $\gamma \subset \P\Omega$ denote a smooth connected arc and $\Sigma$ be the graph of a smooth function $u(x)$ on $\Omega$ in the product manifold $N\PLH\R$.  \\
\indent  Suppose $\Sigma$ is complete approaching to $\gamma$. Then we have 
\begin{enumerate} 
	\item if $\Sigma$ is a translating or minimal graph, then $\gamma$ is a geodesic arc;
	\item if $\Sigma$ is a CMC graph, then $\gamma$ has constant principle curvature. 
\end{enumerate} Moreover only one of the following holds: (1) $u(x)\rightarrow +\infty$ as $x\rightarrow x_0$ for all $x_0\in \gamma$; (2) $u(x)\rightarrow -\infty$ as $x\rightarrow x_0$ for all $x_0\in \gamma$. 
\et 
\br In the case that $N^2$ is $\R^2$ and $\Sigma$ is a complete translating graph in $\R^3$, the above result is obtained by Shahriyari \cite{Sha15}. 
\er 
In the sequel we only prove Theorem \ref{thm:geo} in the case of translating graphs. The proof in the cases of minimal graphs and minimal graphs only requires some minor modifications (see Remark \ref{rm:ce} and Remark \ref{rm:last}).\\
\indent 
From now on, we denote the graph of $u(x)$ by $\Sigma$ which is a translating graph over $\Omega$. Fix a point $x_0\in \gamma$ and let $U_{x_0}$ be an open bounded neighborhood of $x_0$ in $N^2$ such that its intersection with $\gamma$ is a connected arc passing through $x_0$. This arc is written as $\gamma_{x_0}$.  \\
\indent To show Theorem \ref{thm:geo}, our objective is to show that $\gamma_{x_0}$ is geodesic and $\{u(x_n)\}$ has the property as described in the theorem when $\{x_n\}$ approaches to the points on $\gamma_{x_0}$. \\
\indent We start with the following result based on the curvature estimate in the previous section. Notice that $U_{x_0}$ has a compact closure in $N^2$. 
\bl\label{lm:ti} For any point $p=(y,u(y))$ on the translating graph $\Sigma$ where $y\in U_{x_0}$, then $\Sigma$ is a graph (in exponential coordinate of $N^2\PLH\R$ ) over the disk $D_{\delta}\subset T_p\Sigma$ of radius $\delta$. Such graph is denoted by $G(p)$. Moreover $\delta$ and the geometry of $G(p)$ only depend on $U_{x_0}$. The conclusion is also valid when $\Sigma$ is a minimal graph or CMC graph in $N^2\PLH\R$. 
\el 
\br\label{rm:ce} This is the only one place we apply the curvature estimate of translating graphs in Thoerem \ref{thm:est2}. The conclusions in this lemma for minimal graphs and CMC graphs follow from  Theorem \ref{thm:est} and Theorem \ref{thm:cmc:ce} respectively with a similar derivation. Notice that a CMC graph is stable by Zhang \cite{Zhang05}. A similar application of Zhang's estimate can be found in Hauswirth-Rosenberg-Spruck \cite{HRS08}. 
\er 
\bp Since $U_{x_0}$ has a compact closure, then the sectional curvature of $N^2$ on $U_{x_0}$ satisfies  $|K_N |\leq k^2_{x_0}-\F{1}{4}$ for some constant $k_{x_0}$ only depending on $U_{x_0}$.\\
\indent For any $p\in \Sigma$, the exponential map 
$$
\exp_{p}: B_{r_1}(0)\rightarrow \NR
$$
will be a diffeomorphism on the ball in $T_{p}(\NR)=\R^3$ centered at $0$ with radius $r_1$. Here this $r_1$ only depends on $U_{x_0}$ and is independent of $p$. We can equip a metric in $B_{r_1(0)}$ such that the exponential map is an local isometry and $\Sigma$ is a graph near the origin over $T_p\Sigma \cap B_{r_1(0)}(p)$. \\
\indent On the other hand, according to Theorem \ref{thm:est2}, there is a ball in $\Sigma$ with radius $r_0e^{-1}-r$ such that the second fundamental form of $\Sigma$ is uniformly bounded above by $Cr^{-2}$. Let $r=\F{1}{2}r_0e^{-1}$ and $\delta=\F{1}{2}r_0e^{-1}$. Since the exponential map is a uniformly local geometry, we obtain the disk $D_\delta$ in $B_{r_1}(0)$. The geometry of $G(p)$ is determined by the second fundamental form which also only depends on $U_{x_0}$. 
\ep 
For each point $p\in \Sigma$, we translate vertically the graph $G(p)$ into the slice $N^2\PLH\{0\}$ as follows. Let $p=(x^*,u(x^*)$ and $G(p)$ be given in Lemma \ref{lm:ti} with a  representation $(x,u(x))$ where $x$ belongs to some open set $U$. Then its \textbf{vertically translating graph} $F(p)$ is given by 
$$
(x,u(x)-u(x^*)) \subset \NR \text{\quad where}\quad (x,u(x))\in G(p)
$$
This operation does not change any geometric property of $G(p)$.\\
\indent  For any sequence $\{x_n\}\in \Omega$ converging to $x_0$, we conclude that the sequence of $\{u(x_n)\}$ is unbounded. Otherwise the completeness of $\Omega$ approaching to $\gamma$ implies that $x_0\in \Omega$. It is a contradiction.\\
\indent  Let $\{p_n\}$ be the sequence $\{(x_n,u(x_n))\}$ on $\Sigma$ where $\{x_n\}$ converges to $x_0$ as $n\rightarrow \infty$. 
Let $F(p_n)$ and $G(p_n)$ be defined as above with the $\delta$ given in Lemma \ref{lm:ti}. \\
\indent 
The case of translating graphs in Theorem \ref{thm:geo} can be concluded from  the following lemma and the connectedness of $\gamma$. 
\bl \label{lm:minimal} Let $\Sigma$ be a translating graph given in Lemma \ref{lm:ti}. After choosing subsequence, the sequence $F(p_n)$ converges uniformly to $\Gamma\PLH [-\F{\delta}{2},\F{\delta}{2}]\subset \NR$ in the $C^2$ topology where  $\Gamma\subset \gamma_{x_0}$ is a connected geodesic for sufficiently small $\delta$. Moreover only one of the followings holds:
\begin{enumerate} 
	\item  $u(x_n)\rightarrow +\infty$ as $x_n\rightarrow x$ for all $x\in \Gamma$;
	\item $u(x_n)\rightarrow -\infty$ as $x_n\rightarrow x$ for all $x\in \Gamma$.
\end{enumerate}
\el    
\bp We can choose $\delta$ sufficiently small if necessary. By Lemma \ref{lm:ti}, $F(p_n)$, translating $G(p_n)$ into the slice $N^2\PLH\{0\}$, have bounded uniform geometry. After choosing a subsequence, $F(p_n)$ will converge uniformly to a connected surface $F$ passing $x_{0}$ in the $C^2$ topology by Theorem \ref{thm:est2}.\\
\indent Now we claim that the normal vector of $F$ at  $x_0$ is orthogonal to $\P_{r}$, i.e $\Ta =0$. Otherwise the angle function $\Ta$ on $F_n$ shall have positive lower bound. Notice that the diameter of $F_n $ is  $\delta$. When $x_n$ is sufficiently close to $x_0$, $F_n$ has to contain some point in $N^2\PLH\R$ such that its projection into $N^2$ lies outside $\Omega$ because of the two facts mentioned above. It contradicts to the fact that $\Sigma$ is complete approaching to $\gamma$. \\
\indent  By Theorem \ref{thm:convergence} we have $\Ta\equiv 0$ on $F$. Because $H=\Ta$, $F$ is also minimal in $N^2\PLH\R$.  We denote the intersection between $F$ and the slice $N^2\PLH\{0\}$ by $\Gamma$. Thus $\Gamma$ is a geodesic in $N^2$ and $F=\Gamma\PLH [-\F{\sigma}{2},\F{\sigma}{2}]$. Moreover $\Gamma$ is connected because $F$ is connected. \\
\indent  Let $\Gamma_n$ be the intersection between $F_n$ and $N^2\PLH\{0\}$ belonging to $\bar{\Omega}$. According to the definition of $F_n$, we have 
$$
\Gamma_n=\{(x,u(x))\in G_n, u(x)=u(x_n) \}
$$ 
Since $F_n$ converges to $F$, $\Gamma_n$ converges to $\Gamma$ as $n\rightarrow \infty$. We conclude that $\Gamma\subset \bar{\Omega}$.\\
\indent Suppose $\Gamma$ does not belongs to $\gamma_x$. Then there is a sequence $\{y_n\in\Gamma_n\}$ converges to $y\in \Omega$ on $\Gamma$. Thus the sequence $\{u(y_n)\}$ converges to $u(y)$ which is a finite number. This implies that the sequence $\{u(y_n)=u(x_n)\}$ is bounded from the definition of $\Gamma_n$. This is a contradiction to the fact that $\{u(x_n)\}$ is an unbounded sequence. As a result $\Gamma\subset \gamma_x$. \\
\indent Assuming that $u(x_n)\rightarrow +\infty$ as $x_n \rightarrow x_0$. Then  for any sequence $\{x_n'\}$ approaching to $x_0$, then $u(x'_n)\rightarrow +\infty$. Otherwise by the intermediate theorem of continuous functions, there is a sequence $\{x_n^{''}\}$ such that $u(x_n^{''})$ converges to a finite number as $x_n^{''}\rightarrow x_0$. Again it contradicts the completeness of $\Sigma$ approaching to $\gamma$. The word-by-word derivation also works for any point $y\in \Gamma$. Thus we conclude (1). The proof of (2) is similar when we assume $u(x_n)\rightarrow-\infty$ as $x_n\rightarrow x_0$. The proof of Lemma \ref{lm:minimal} is complete.  \ep 
\br\label{rm:last} The conclusion in Lemma \ref{lm:minimal} is also valid when $\Sigma$ is minimal or CMC. The only modification is that $\Gamma$ is a geodesic or an arc  with constant principle curvature respectively. 
\er


\section{Acknowledgement}{This work was supported by the National Natural Science Foundation of China, No.11261378 and No.11521101. The author is very grateful to the encouragement from Prof. Lixin Liu. The author also thanks the referees for careful readings and helpful suggestions.}





\begin{appendix}
	\section{Examples of translating graphs}
	In this section we construct some examples of translating graphs to mean curvature flows when the surface $N^2$ has a domain with certain special warped product structure. \\
	\indent  Suppose $N^2$ is a complete Riemannian surface with a metric $\sigma$ containing a domain $N^2_0$ equipped with the following coordinate system:
	\be \label{metric:structure}
	\{\theta \in S^1, r\in [0, r_0)\}\quad\text{with}\quad \sigma=dr^2+h^2(r)d\ta^2
	\ene
	where $d\ta^2$ is the standard metric on the unit circle $S^1$, $h(r)$ is a positive function satisfying $h(0)=0$, $h'(0)=1$ with $h'(r)\neq 0$ for all $r\in (0, r_0)$. For more detail on warped product metric, we refer to Section 2 in \cite{HZ16}.\\
	\indent The following result discusses the existence of translating graphs in $\NR$ with the structure in \eqref{metric:structure}.
	\bt Let $N^2$ be a surface mentioned above. Let $u(r):[0,r_0)\rightarrow \R$ be a $C^2$ solution of the following ordinary equation
	\be\label{eq:u}
	\F{u_{rr}}{1+u_r^2}+\F{h'(r)}{h(r)}u_r=1
	\ene
	with $u_r(0)=0$ for $r\in [0,r_0)$. Then $\Sigma=(x,u(r))$ for $r\in [0,r_0)$ is a translating graph in $\NR$ where $x=(r,\ta)\in N^2_0$ given by \eqref{metric:structure}. If $r_0=\infty$, then $\Sigma$ is complete.
	\et
	\br \eqref{eq:u} is an ODE. The existence of its solution is obvious.\\
	\indent Two concrete examples are the unit sphere $S^2$ and the hyperbolic plane $\H^2$. In the former case, $N^2_0$ is the hemisphere with $h(r)=\sin(r)$ for $r\in [0,\F{\pi}{2})$ and the spherical metric is written as 
	$
	dr^2+\sin(r)d\ta^2 
	$. In the case of $\H^2$ the hyperbolic metric is written as 
	$
	dr^2+\sinh(r)d\ta^2
	$. 
	\er 
	\bp Suppose $r_0=\infty$. Then $N^2_0$ is simply connected and should be a whole $N^2$. Thus $\Sigma$ is complete.\\
	\indent Now we show that $\Sigma$ is a translating graph.\\
	\indent According to \eqref{eq:angle} it is sufficient to derive the identity 
	\be \label{def:TS}
	H=-\Ta
	\ene  
	where $H$ is the mean curvature of $\Sigma$ and $\vec{v}$ is its upward normal vector.\\
	\indent Fix a point $(x,u(x))$ on $\Sigma$ where $x\in N_0^2$ and the polar coordinate of $x$ in $N_0^2$ is not $(0,0)$. There is a natural frame $\{\P_r,\P_\ta\}$ according to the polar coordinate on $N^2_0$ by \eqref{metric:structure}. Let $\Sigma$ be the graph of $u(x)=u(r)$ in $\NR$. Let $u_r, u_\ta $ denote the partial derivative of $u$. Thus we have a natural frame $\{X_1=\P_r+u_r\P_3,X_2=\P_\ta\}$ on $\Sigma$.  Here $\P_3$ denotes the vector field tangent to $\R$. We also use the fact $u_\ta=0$. Then the metric on $\Sigma$ and the upward normal vector of $\Sigma$ are given by
	\begin{align*}
	g_{11}&=\la X_1, X_1\ra=1+u_r^2,\quad  g_{12}=\la X_1, X_2\ra =0 \\
	g_{22}&=\la X_2, X_2\ra=h^2(r)\\
	\vec{v}& =\F{\P_3-u_r\P_r}{\sqrt{1+u_r^2} }                                                                 
	\end{align*}
	Let $\bnb$ denote the covariant derivative of $\NR$. Then its second fundamental form is 
	$$
	h_{11}=-\la\bnb_{X_1}X_1,\vec{v}\ra =\F{u_{rr}}{\sqrt{1+u_r^2}}\quad
	h_{22}=\la \bnb_{\P_\ta}\P_\ta,\vec{v}\ra =-h'(r)h(r)\F{u_r}{\sqrt{1+u_r^2}}
	$$
	where we use the fact $\la \bnb_{\P_\ta}\P_\ta,\P_r\ra=-h'(r)h(r)$ (for more detail see Section 2 in \cite{HZ16}). Then the mean curvature of $\Sigma$ with respect to $\vec{v}$ 
	is 
	$$
	H=g^{11}h_{11}+g^{22}h_{22}=-\F{1}{\sqrt{1+u_r^2}}(\F{u_{rr}}{1+u_r^2}+\F{h'(r)}{h(r)}u_r)=-\F{1}{\sqrt{1+u_r^2}}
	$$
	by \eqref{eq:u}. On the other hand we have  
	$$
	\Ta =\la\vec{v},\P_3\ra
	=\F{1}{\sqrt{1+u_r^2}}=-H
	$$ Hence $\Sigma$ is a translating graph. The proof is complete. 
	\ep

\end{appendix}
\bibliographystyle{abbrv}	
\bibliography{Ref_Thesis}

\begin{thebibliography}{10}

\bibitem{AW94}
S.~J. Altschuler and L.~F. Wu.
\newblock Translating surfaces of the non-parametric mean curvature flow with
  prescribed contact angle.
\newblock {\em Calculus of Variations and Partial Differential Equations},
  2(1):101--111, 1994.

\bibitem{AngV97}
S.~B. Angenent and J.~J.~L. Velazquez.
\newblock Degenerate neckpinches in mean curvature flow.
\newblock {\em J.reine Angew.math}, 1997(482):15--66, 1997.

\bibitem{Ang95}
S.~B. Angenent and J.~J.~L. Vel�zquez.
\newblock Asymptotic shape of cusp singularities in curve shortening.
\newblock {\em Duke Mathematical Journal}, 77(1):71--110, 1995.

\bibitem{CSS07}
J.~Clutterbuck, O.~C. Schn\"urer, and F.~Schulze.
\newblock Stability of translating solutions to mean curvature flow.
\newblock {\em Calc. Var. Partial Differential Equations}, 29(3):281--293,
  2007.

\bibitem{CM02}
T.~H. Colding and W.~P. Minicozzi.
\newblock Estimates for parametric elliptic integrands.
\newblock {\em International Mathematics Research Notices}, 2002(6):291--297,
  2002.

\bibitem{CR10}
P.~Collin and H.~Rosenberg.
\newblock Construction of harmonic diffeomorphisms and minimal graphs.
\newblock {\em Ann. of Math. (2)}, 172(3):1879--1906, 2010.

\bibitem{doC92}
M.~P.~a. do~Carmo.
\newblock {\em Riemannian geometry}.
\newblock Mathematics: Theory \& Applications. Birkh\"auser Boston, Inc.,
  Boston, MA, 1992.
\newblock Translated from the second Portuguese edition by Francis Flaherty.

\bibitem{Ecm10}
M.~Eichmair.
\newblock The plateau problem for marginally outer trapped surfaces.
\newblock {\em Journal of Differential Geometry}, 83(3):551--584, 2010.

\bibitem{EM16}
M.~Eichmair, J.~Metzger, et~al.
\newblock Jenkins-{S}errin-type results for the {J}ang equation.
\newblock {\em Journal of Differential Geometry}, 102(2):207--242, 2016.

\bibitem{GJJ10}
C.~Gui, H.~Jian, and H.~Ju.
\newblock Properties of translating solutions to mean curvature flow.
\newblock {\em Discrete Contin. Dyn. Syst.}, 28(2):441--453, 2010.

\bibitem{HRS08}
L.~Hauswirth, H.~Rosenberg, and J.~Spruck.
\newblock On complete mean curvature {$\frac 12$} surfaces in {$\Bbb
  H^2\times\Bbb R$}.
\newblock {\em Comm. Anal. Geom.}, 16(5):989--1005, 2008.

\bibitem{HP99}
G.~Huisken and A.~Polden.
\newblock Geometric evolution equations for hypersurfaces.
\newblock {\em Lecture Notes in Mathematics}, 1713:45--84, 1999.

\bibitem{HS08}
G.~{Huisken} and C.~{Sinestrari}.
\newblock {Mean curvature flow with surgeries of two-convex hypersurfaces}.
\newblock {\em Inventiones Mathematicae}, 175:137--221, Sept. 2008.

\bibitem{JS68}
H.~Jenkins and J.~Serrin.
\newblock The {D}irichlet problem for the minimal surface equation in higher
  dimensions.
\newblock {\em J.reine Angew.math}, 1968(229):170--187, 1968.

\bibitem{PL09}
A.~L. Pinheiro.
\newblock A {J}enkins-{S}errin theorem in {$M^2\times\Bbb R$}.
\newblock {\em Bull. Braz. Math. Soc. (N.S.)}, 40(1):117--148, 2009.

\bibitem{Soe83}
R.~Schoen.
\newblock Estimates for stable minimal surfaces in three dimensional manifolds.
\newblock In {\em Seminar on minimal submanifolds}, volume 103, pages 111--126.
  Princeton Univ. Press. Princeton, 1983.

\bibitem{Sha15}
L.~Shahriyari.
\newblock Translating graphs by mean curvature flow.
\newblock {\em Geometriae Dedicata}, 175(1):57--64, 2015.

\bibitem{spr72}
J.~Spruck.
\newblock Infinite boundary value problems for surfaces of constant mean
  curvature.
\newblock {\em Arch. Rational Mech. Anal.}, 49:1--31, 1972/73.

\bibitem{Sj16}
J.~{Sun}.
\newblock {Lagrangian $L$-stability of Lagrangian Translating Solitons}, Dec.
  2016.
\newblock arXiv:1612.06815.

\bibitem{WXJ11}
X.-J. Wang.
\newblock Convex solutions to the mean curvature flow.
\newblock {\em Ann. of Math. (2)}, 173(3):1185--1239, 2011.

\bibitem{Zhang05}
S.~Zhang.
\newblock Curvature estimates for cmc surfaces in three dimensional manifolds.
\newblock {\em Mathematische Zeitschrift}, 249(3):613--624, 2005.

\bibitem{HZ16}
H.~Zhou.
\newblock Inverse mean curvature flows in warped product manifolds.
\newblock {\em The Journal of Geometric Analysis}, Jun 2017.

\end{thebibliography}

\end{document}